\documentclass[11pt,leqno]{article} 
\usepackage{graphics}
\newtheorem{thm}{Theorem}[section]
\newtheorem{lma}{Lemma}[section]

\newcommand{\beqa}{\begin{eqnarray}}
\newcommand{\eeqa}{\end{eqnarray}}

\newcommand{\pf}{\noindent {\bf Proof:} $\s$ }
\newcommand{\epf}{ \hfill$\diamondsuit$ \medskip}

\newcommand{\beq}{\begin{equation}}
\newcommand{\eeq}{\end{equation}}
\newcommand{\lbl}{\label}
\newcommand{\s}{\; \;}

\newcommand{\ep}{\epsilon}

\newcommand{\la}{\lambda}
\newcommand{\mb}{\mbox}
\newcommand{\ra}{\rightarrow}
\newcommand{\al}{\alpha}
\newcommand{\om}{\omega}

\title{Global solution curves for self-similar equations}

\author{
Philip Korman   \\ 
Department of Mathematical Sciences \\ 
University of Cincinnati \\ 
Cincinnati Ohio 45221-0025 \\
}

\date{}

\begin{document}

\maketitle

\begin{abstract} 
\noindent
We consider positive solutions of a semilinear Dirichlet problem
\[
\Delta u+\la f(u)=0, \s \mbox{for $|x|<1$}, \s u=0 , \s \mbox{when $|x|=1$}
\]
on a unit ball in $R^n$.  For four classes of self-similar equations it is possible to parameterize the entire (global) solution curve through the solution of a single initial value problem. This allows us to derive results on the multiplicity of solutions, and on their Morse indices. In particular, we easily recover the classical results of D.D. Joseph and T.S. Lundgren \cite{JL} on the Gelfand problem. Surprisingly, the situation turns out to be different for the generalized Gelfand problem, where infinitely many turns are possible for any space dimension $n \geq 3$. We also derive detailed results for the equation modeling electrostatic micro-electromechanical systems (MEMS), in particular we easily recover  the main result of Z. Guo and J. Wei \cite{GW}, and we show that the Morse index of the solutions increases by one at each turn. We also consider the self-similar Henon's equation.
 \end{abstract}

\begin{flushleft}
\end{flushleft}

\begin{flushleft}
\end{flushleft}

\section{Introduction}

We consider radial solutions on a ball in $R^n$ for four special classes of equations, the ones self-similar under scaling.
For example, consider the so called Gelfand equation ($u=u(x)$, $x \in R^n$)
\beq
\lbl{0.1}
\Delta u+\la e^u=0, \s \mbox{for $|x|<1$}, \s u=0 , \s \mbox{when $|x|=1$} \,.
\eeq
Here $\la$ is a positive parameter. By the maximum principle, solutions of (\ref{0.1}) are positive, and then by the classical theorem of B. Gidas, W.-M. Ni and L. Nirenberg \cite{GNN} they are radially symmetric, i.e., $u=u(r)$, $r=|x|$, and it satisfies
\beq
\lbl{0.2}
u''+\frac{n-1}{r}u' + \la e^u=0, \s \mb{for $0<r<1$}, \s u'(0)=u(1)=0 \,.
\eeq
This theorem also asserts that $u'(r)<0$ for all $0<r<1$, which implies that the value of $u(0)$ gives the $L^{\infty}$ norm of our solution. Moreover, $u(0)$ is a {\em global parameter}, i.e., it uniquely identifies the solution pair $(\la ,u(r))$, see e.g., P. Korman \cite{K1}.
It follows that a two-dimensional curve $(\la,u(0))$ completely describes the solution set of (\ref{0.1}). The change of variables $v= u+a$, $\xi= br$, with constant $a$ and $b$ will transform the equation in (\ref{0.2}) into the same equation if $e^a=b^2$. Here is what this self-similarity ``buys" us.
Let $w(t)$ be the solution of the following initial value problem
\beq
\lbl{0.21}
w''+\frac{n-1}{t}w'+  e^{w}=0,   \s  \s w(0)=0, \s w'(0)=0 \s (t>0)\,,
\eeq
which is easily seen to be negative, and defined for all $t \in (0,\infty)$. It turns out that $w(t)$ gives us the entire solution curve of (\ref{0.2}) (or of (\ref{0.1})):
\beq
\lbl{0.3}
\left(\la,u(0) \right)=\left(t^2e^{w(t)}\,, -w(t) \right) \,,
\eeq
parameterized by $t \in (0,\infty)$. In particular, $\la=\la(t)=t^2e^{w(t)}$, and 
\[
\la '(t)=t e^{w} \left(2+tw' \right) \,,
\]
so that the solution curve travels to the right (left) in the $(\la,u(0))$ plane if  $2+tw'>0$ ($<0$). This makes us  interested in the roots of the function $2+tw'$. If we set this function to zero
\[
2+tw'=0 \,,
\]
then solution of this equation is of course $w(t)=a-2 \ln t$. Amazingly, if we choose $a=\ln (2n-4)$, $n \geq 3$, then 
\[
w_0(t)=\ln (2n-4)-2 \ln t
\]
is also a solution of the equation in (\ref{0.21})! We show that $w(t)$ tends to $w_0(t)$ as $t \ra \infty$, and the issue turns out to be  how many times $w(t)$ and $w_0(t)$ cross as $t \ra \infty$. We propose to call  $w(t)$ {\em the generating solution}, and $w_0(t)$ {\em the guiding solution}. We show that the solution curve makes infinitely many turns if and only if $w(t)$ and $w_0(t)$ intersect infinitely many times. Then we prove that for $3\leq n \leq 9$, $w(t)$ and $w_0(t)$ intersect infinitely many times, and hence the solution curve makes infinitely many turns, which is a part of the  classical result of  D.D. Joseph and T.S. Lundgren \cite{JL}, see also J. Bebernes and D. Eberly \cite{BE} for an exposition.  D.D. Joseph and T.S. Lundgren \cite{JL} also proved that for $n=1,2$ the solution curve makes exactly one turn, while for $n \geq 10$ there are no turns (we recover this result for $n \geq 10$ too). Our approach provides a remarkably short route to  the classical result of D.D. Joseph and T.S. Lundgren \cite{JL}, and some new results for other equations.
\medskip

A similar approach works for the radially symmetric solutions of the generalized Gelfand's problem
\[
\Delta u+\la |x|^{\al} e^u=0, \s \mbox{for $|x|<1$}, \s u=0 , \s \mbox{when $|x|=1$} \,,
\]
with a constant $\al>0$. Remarkably, the picture here turns out to be different! We show that the solution curve makes infinitely many turns, provided that 
\[
3 \leq n <10+4 \al \,.
\]
We see that unlike the Gelfand equation, infinitely many turns occur for any space dimension $n \geq 3$, for large enough $\al$. This result is sharp, because we prove that there are at most   two turns if  $n \geq 10+4 \al $.

A similar approach works for three  other classes of equations, notably 
\beq
\lbl{0.3a}
\s \s u''+\frac{n-1}{r}u' + \la \, \frac{r^{\alpha}}{(1-u)^p} =0,  \s u'(0)=u(1)=0 \,, \s 0<u(r)<1  \,,
\eeq
 modeling the electrostatic micro-electromechanical systems (MEMS), see J.A. Pelesko \cite{P}, N. Ghoussoub and Y.  Guo \cite{GG}, and  Z. Guo and J. Wei \cite{GW} for some of the active recent research. We give a much shorter and more elementary proof of the main result of Z. Guo and J. Wei \cite{GW}, which states that the solution curve makes infinitely many turns, provided that 
\beq
\lbl{0.5}
2 \leq n <2+ \frac{2(\al +2)}{p+1} \, \left(p+\sqrt{p^2+p} \right)\,.
\eeq
(If $2 \leq n \leq 6$, this inequality holds for all $\al \geq 0$ and $p>1$.)
We also show that outside of this range the solution curve makes at most two  turns, which is a new result.
\medskip

For the Gelfand problem (\ref{0.1}) it was shown by K. Nagasaki and T.  Suzuki \cite{N2} that at each turn of the solution curve the Morse index of the solution increases by one. We recover this result, and then prove that the same thing is true for the MEMS problem (\ref{0.3a}), which is a new result.
Finally, we use the self-similar nature of Henon's equation to discuss the exact multiplicity of the symmetry breaking solutions.

\section{Parameterization of the global solution curves}
\setcounter{equation}{0}
\setcounter{thm}{0}
\setcounter{lma}{0}

Consider the problem 
\beq
\lbl{1}
u''+\frac{n-1}{r}u' + \la \, \frac{r^{\alpha}}{(1-u)^p} =0 \s\s \mb{for $0<r<1$} \,, 
\eeq
\[
u'(0)=u(1)=0 \,, \s 0<u(r)<1  \,,
\]
which arises in modeling of electrostatic micro-electromechanical systems (MEMS), see \cite{P}, \cite{GG}, \cite{GW}. Here $\la $ is a positive parameter, $\al >0$ and $p>1$ are constants. Any  solution $u(r)$ of (\ref{1}) is a positive and  decreasing function (by the maximum principle), so that $u(0)$ gives its maximum value. It is known, see e.g.,  P. Korman \cite{K2}, that $u(0)$ is a {\em global parameter}, i.e., it uniquely identifies the solution pair $(\la ,u(r))$ (the proof is by scaling).
It follows that a two-dimensional curve $(\la,u(0))$ completely describes the solution set of (\ref{1}). Our goal is to compute the global solution curve $(\la,u(0))$. Let $1-u=v$. Then $v(r)$ satisfies
\beq
\lbl{2}
v''+\frac{n-1}{r}v'= \la \, \frac{r^{\al}}{v^p}  \s\s \mb{for $0<r<1$,} \s v'(0)=0, \s v(1)=1 \,.
\eeq
Assume that $v(0)=a$. We scale $v=aw$, and $t=br$. The constants $a$ and $b$ are assumed to satisfy
\beq
\lbl{3}
\la =a^{p+1} b ^{\al +2} \,.
\eeq
Then (\ref{2}) becomes
\beq
\lbl{4}
w''+\frac{n-1}{t}w'=  \frac{t^{\al}}{w^p}, \s w(0)=1, \s w'(0)=0 \,.
\eeq
It is easy to see that the solution $w(t)$ of (\ref{4}) is an increasing function, defined for all $t>0$. We can compute it numerically (over a large interval). It turns out that this  particular solution $w(t)$ gives us the entire solution curve of (\ref{1})! We have
\[
1=v(1)=aw(b) \,,
\]
and so $a=\frac{1}{w(b)}$, and then $\la =\frac{b ^{\al +2}}{w^{p+1}(b)}$. The global solution curve  is
\beq
\lbl{4a}
(\la,u(0))=\left(\frac{b ^{\al +2}}{w^{p+1}(b)}\,, 1-\frac{1}{w(b)} \right) \,,
\eeq
parameterized by $b \in (0,\infty)$. This parameterization was pointed out previously in \cite{P},  and was then  used in  \cite{GG}.
\medskip

\noindent
{\bf Example} Using {\em Mathematica}, we have solved the problem (\ref{1}) with $p=2$, $n=2$ and $\al =0.2$.
The global solution curve, obtained through the parameterization (\ref{4a}), is given in Figure 1.

\begin{figure}
\begin{center}
\scalebox{0.6}{\includegraphics{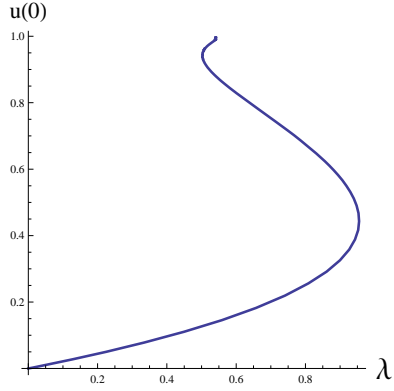}}
\caption{  Solution curve for the problem (\ref{1})}
\end{center}
\end{figure}

\medskip

We now show that  the situation is similar for three other important classes of equations.
Consider the problem 
\beq
\lbl{8}
\s\s\s u''+\frac{n-1}{r}u' + \la \, r^{\alpha}(1+u)^p =0 \s \mb{for $0<r<1$,} \s u'(0)=0, \s u(1)=0 \,.
\eeq
We set $v=1+u$, followed by $v=aw$, and $t=br$, where $a=v(0)=1+u(0)$.  
The constants $a$ and $b$ are assumed to satisfy
\beq
\lbl{9}
\la =\frac{ b ^{\al +2}}{a^{p-1}} \,.
\eeq
Then (\ref{8}) becomes
\beq
\lbl{10}
w''+\frac{n-1}{t}w'+  t^{\al}w^p=0, \s w(0)=1, \s w'(0)=0 \,.
\eeq
The solution of this problem is easily seen to be a decreasing function, which, for sub-critical $p$, vanishes at some $t_0>0$. (If $p \geq \frac{n+2+2\alpha}{n-2}$, then $w(t)$ has no roots on $(0,\infty)$, see e.g., T. Kusano and M. Naito \cite{KN}, or E. Yanagida and S. Yotsutani \cite{Y}.) As before, $a=\frac{1}{w(b)}$, and then $\la =b ^{\al +2}w^{p-1}(b)$.
The global solution curve  is
\[
(\la,u(0))=\left(b ^{\al +2}w^{p-1}(b)\,, -1+\frac{1}{w(b)} \right) \,,
\]
parameterized by $b \in (0,t_0)$.
\medskip

Next, we consider the generalized Gelfand's equation
\beq
\lbl{12}
\s\s\s u''+\frac{n-1}{r}u' + \la \, r^{\alpha}e^u =0 \s \mb{for $0<r<1$,} \s u'(0)=0, \s u(1)=0 \,.
\eeq
We set $u=w+a$, $t=br$, with $a=u(0)$. The constants $a$ and $b$ are assumed to satisfy
\[
\la = b ^{\al +2}e^{-a} \,.
\]
Then (\ref{12}) becomes
\[
w''+\frac{n-1}{t}w'+  t^{\al}e^w=0, \s w(0)=0, \s w'(0)=0 \,.
\]
We compute numerically the solution of this problem $w(t)$, which  is a negative decreasing function, defined for all $t>0$. We have
\[
0=u(1)=a+w(b) \,,
\]
i.e., $a=-w(b)$, and then  $\la = b ^{\al +2}e^{w(b)}$. The global solution curve for (\ref{12}) is 
\[
(\la,u(0))=\left(b ^{\al +2}e^{w(b)}\,, -w(b) \right),
\]
parameterized by $b \in (0,\infty)$. Moreover, $u(r)=-w(b)+w(br)$ is the solution of (\ref{12}) at $\la = b ^{\al +2}e^{w(b)}$.
\medskip

Finally, we consider 

\beq
\lbl{14}
\s\s\s u''+\frac{n-1}{r}u' + \la \, r^{\alpha}e^{-u} =0 \s \mb{for $0<r<1$,} \s u'(0)=0, \s u(1)=0 \,.
\eeq
We set $u=w+a$, $t=br$, with $a=u(0)$. The constants $a$ and $b$ are assumed to satisfy
\[
\la = b ^{\al +2}e^{a} \,.
\]
Then (\ref{14}) becomes
\[
w''+\frac{n-1}{t}w'+  t^{\al}e^{-w}=0, \s w(0)=0, \s w'(0)=0 \,.
\]
We compute numerically the solution of this problem $w(t)$, which  is a negative decreasing function, tending to $-\infty$ at some   $t_1>0$. We have
\[
0=u(1)=a+w(b) \,,
\]
i.e., $a=-w(b)$, and then  $\la = b ^{\al +2}e^{-w(b)}$. The global solution curve for (\ref{14}) is 
\[
(\la,u(0))=\left(b ^{\al +2}e^{-w(b)}\,, -w(b) \right),
\]
parameterized by $b \in (0,t_1)$. 

\section{A generalization of Joseph and Lundgren's result}
\setcounter{equation}{0}
\setcounter{thm}{0}
\setcounter{lma}{0}

As we saw above, for  Gelfand's problem 
\beq
\lbl{16}
\s\s\s u''+\frac{n-1}{r}u' + \la \, e^u =0 \s \mb{for $0<r<1$,} \s u'(0)=0, \s u(1)=0 
\eeq
the solution curve $(\la, u(0))$ is given by $\left(t ^{2}e^{w(t)}\,, -w(t) \right)$, parameterized by $t \in (0,\infty)$, where $w(t)$ is the solution of
\beq
\lbl{17}
w''+\frac{n-1}{t}w'+  e^{w}=0,   \s  \s w(0)=0, \s w'(0)=0 \s (t>0)\,.
\eeq
In particular, $\la =t ^{2}e^{w(t)}$, and the issue is how many times this function changes its direction of monotonicity for $t \in (0,\infty)$.
Compute
\beq
\lbl{17.1}
\la '(t)=t e^{w} \left(2+tw' \right),
\eeq
so that we are interested in the roots of the function $2+tw'$. If we set this function to zero
\[
2+tw'=0 \,,
\]
then solution of this equation is of course $w(t)=a-2 \ln t$. Amazingly, if we choose $a=\ln (2n-4)$, $n \geq 3$, then 
\[
w_0(t)=\ln (2n-4)-2 \ln t
\]
is a solution of the equation in (\ref{17})! It turns out that the solution of (\ref{17}) tends to $w_0(t)$ as $t \ra \infty$, and the issue is how many times $w(t)$ and $w_0(t)$ cross as $t \ra \infty$. We shall only consider $n \geq 3$, since for $n=1,2$ the problem (\ref{16}) can be explicitly solved, see e.g., \cite{BE}.

\begin{lma}\lbl{lma:3}
Assume that $w(t)$ and $w_0(t)$ intersect infinitely many times. Then the solution curve of (\ref{16}) makes infinitely many turns.
On the other hand, assume that $w(t)$ and $w_0(t)$ intersect  only a finite number of times, and $(w-w_0)'(t)$ is of one sign for $t>T$, with a point $T$ coming after the last point of intersection. Then  $\la (t)$ is monotone for $t>T$.
\end{lma}

\pf
Indeed, assuming that $w(t)$ and $w_0(t)$ intersect infinitely many times, let $\{ t_n \}$ denote the points of intersection. At $\{ t_n \}$'s, $w(t)$ and $w_0(t)$ have different slopes (by uniqueness for initial value problems). Since $2+t_nw'_0(t_n)=0$, it follows that $2+t_nw'(t_n)>0$ ($<0$) if  $w(t)$ intersects $w_0(t)$ from below (above) at $t_n$. Hence, on any interval $(t_n,t_{n+1})$ there is a point $t_0$, where $2+t_0w'(t _0)=0$, i.e., $\la '(t_0)=0$, and $t_0$ is a critical point. Since $\la '(t_n)$ and $\la '(t_{n+1})$ have different signs, the solution curve changes its direction over $(t_n,t_{n+1})$.  
\medskip

On the other hand, assume that $w(t)$ and $w_0(t)$ intersect  only a finite number of times, and $(w-w_0)'(t)$ is of one sign for $t>T$, say $w'(t)>w_0'(t)$. Then $2+tw'(t)>2+tw'_0(t)=0$, i.e., $\la' (t)>0$, and the solution curve is monotone for $t>T$.
\epf

\medskip

The linearized equation for (\ref{17}) is
\[
z''+\frac{n-1}{t}z'+  e^{w}z=0 \,.
\]
At the solution $w=w_0(t)$, this becomes
\beq
\lbl{19}
z''+\frac{n-1}{t}z'+  \frac{2n-4}{t^2}z=0 \,,
\eeq
which is Euler's equation! Its characteristic equation has the roots
\beq
\lbl{19.1}
r=\frac{-n+2 \pm \sqrt{(n-2)(n-10)}}{2} \,.
\eeq
When $3 \leq n \leq 9$, the roots are complex, and hence $z(t)$ changes sign infinitely many times. We shall show  that  $w(t)$ tends to $w_0(t)$, and oscillates infinitely many times around $w_0(t)$, which implies infinitely many turns of the solution curve. For other $n$, the  solution curve turns  at most once.  We obtain a remarkably short route to  the classical result of D.D. Joseph and T.S. Lundgren \cite{JL}.
\medskip

We shall present the details for the more general problem ($\al >0$)
\beq
\lbl{30}
\s\s\s u''+\frac{n-1}{r}u' + \la \, r^{\al} e^u =0 \s \mb{for $0<r<1$,} \s u'(0)=0, \s u(1)=0 \,.
\eeq
The solution curve $(\la, u(0))$ is now given by $\left(t ^{2+\al}e^{w(t)}\,, -w(t) \right)$, parameterized by $t \in (0,\infty)$, where $w(t)$ is the solution of
\beq
\lbl{31}
w''+\frac{n-1}{t}w'+t ^{\al} e^{w}=0,   \s  \s w(0)=0, \s w'(0)=0 \s (t>0)\,.
\eeq
In particular, $\la =t ^{2+\al}e^{w(t)}$, and the issue is how many times this function changes the direction of monotonicity for $t \in (0,\infty)$.
Compute
\[
\la '(t)=t^{\al +1}  e^{w} \left(2+\al+tw' \right) \,,
\]
so that we are interested in the roots of the function $2+\al+tw'$. If we set this function to zero
\[
2+\al+tw'=0 \,,
\]
then the general solution of this equation is $w(t)=a-(2+\al) \ln t$. If we choose $a=\ln (2+\al)(n-2)$, $n \geq 3$, then 
\[
w_0(t)=\ln (2+\al)(n-2)-(2+\al) \ln t
\]
is a solution of the equation in (\ref{31}). It turns out that the solution $w(t)$ of (\ref{31}) tends to $w_0(t)$ as $t \ra \infty$, and the issue is how many times $w(t)$ and $w_0(t)$ cross as $t \ra \infty$.
\medskip

The linearized equation for (\ref{31}) is
\[
z''+\frac{n-1}{t}z'+  t ^{\al}e^{w}z=0 \,.
\]
At the solution $w=w_0(t)$, this becomes
\beq
\lbl{32}
z''+\frac{n-1}{t}z'+  \frac{(2+\al)(n-2)}{t^2}z=0 \,,
\eeq
which is Euler's equation. Its characteristic equation has the roots
\[
r=\frac{-n+2 \pm \sqrt{(n-2)(n-10-4 \al)}}{2} \,.
\]
When $3 \leq n <10+4 \al$, the roots are complex. Hence, $z(t)$ changes sign  infinitely many times. We shall show that this  implies infinitely many turns of the solution curve, while  for other $n$ at most two turns of the solution curve is possible.
\medskip

We shall need the following version of Sturm's comparison theorem.
\begin{lma}\lbl{lma:1}
Consider the following two equations
\beq
\lbl{33}
y''+\frac{n-1}{t}y'+\frac{a(t)}{t^2}y=0 \,,
\eeq
\beq
\lbl{34}
v''+\frac{n-1}{t}v'+\frac{b(t)}{t^2}v=0 \,.
\eeq
Assume that $b(t)>a(t)$ for all $t \in R$. Then $v(t)$ has a root between any two consecutive roots of $y(t)$.
\end{lma}

\pf
Assume that $y(t_1)=y(t_2)=0$, $y(t)>0$ on $(t_1,t_2)$, while on the contrary $v(t)>0$ on $(t_1,t_2)$. From the equations (\ref{33}) and (\ref{34})
\[
\left[t^{n-1} \left(y'v-yv' \right) \right]'=t^{n-3} \left(b(t)-a(t) \right) y(t)v(t) >0 \s \mbox{on $(t_1,t_2)$} \,.
\]
Integrating over $(t_1,t_2)$,
\[
t_2^{n-1} y'(t_2)v(t_2)-t_1^{n-1} y'(t_1)v(t_1)>0 \,,
\]
which is a contradiction, since both terms on the left are non-positive.
\epf

\begin{lma}\lbl{lma:2}
Consider the equation (here $a$, $a_1$, $a_2$ are constants)
\beq
\lbl{35}
y''+\frac{n-1}{t}y'+\frac{a+f(t)}{t^2}y=0 \,.
\eeq
Assume that the equation (\ref{33}), with $a(t)=a$, has infinitely many roots for all $a \in [a_1,a_2]$, while $f(t) \ra 0$, as $t \ra \infty$.
Then the equation (\ref{35}) has infinitely many roots for all $a \in (a_1,a_2)$.
\end{lma}

\pf
Choose an $\ep >0$, so that $a-\ep >a_1$. Since $a+f(t)>a-\ep>a_1$ for $t$ large, the proof follows by Lemma \ref{lma:1}.
\epf

\begin{lma}\lbl{lma:5}
Consider the equation (\ref{35}), with $n \geq 3$, $a>0$, and $\lim _{t \ra \infty} f(t)=0$. Assume that its solution $y(t)$ is bounded on some interval $(t_0,\infty)$. Then $\lim _{t \ra \infty} y(t)=0$.
\end{lma}

\pf
Letting $t=e^s$, we transform (\ref{35})  to
\beq
\lbl{35.1}
y_{ss}+(n-2)y_s+ay=g(s) \,,
\eeq
with $g(s) \equiv -f(e^s)y \ra 0$, as $s \ra \infty$. The roots of the corresponding homogeneous equation are either $\al \pm i \beta$, with $\al <0$, or both roots are negative. Let us assume it is the former case, and the other case is similar. The general solution of (\ref{35.1}) is
\[
y(s)=c_1 e^{\al s} \cos \beta s +c_2 e^{\al s} \sin \beta s+\frac{1}{\beta} \int_0^s e^{\al (s-\xi)} \sin \beta (s-\xi) g(\xi) \, d \xi \,,
\]
and 
\[
|\int_0^s e^{\al (s-\xi)} \sin \beta (s-\xi) g(\xi) \, d \xi| \leq e^{\al s} \int_0^s e^{-\al \xi}  |g(\xi)| \, d \xi \ra 0, \s \mbox{as $s \ra \infty$} \,,
\]
concluding the proof.
\epf

\begin{thm}\lbl{thm:1}
Assume that 
\[
3 \leq n <10+4 \al \,.
\]
Then the solution curve of (\ref{30}) makes infinitely many turns.
\end{thm}

\pf
In view of the Lemma \ref{lma:3}, we need to show that $w(t)$ oscillates infinitely many times around $w_0(t)$. We claim that these functions get close to each other, as $t$ increases. Denote $p(t)=w_0(t)-w(t)$. It satisfies
\beq
\lbl{36}
p''+\frac{n-1}{t}p'+t^{\alpha } a(t) p=0 \,,
\eeq
where $a(t)=\int_0^1  e^{s w_0(t)+(1-s) w(t)} \, ds>0$.
We have $p(\ep)>0$ and $p'(\ep)<0$, for $\ep>0$ small. From (\ref{36}), 
\[
\left( t^{n-1} p' \right)'=-t^{n+\alpha-1} a(t)p<0 \,, \s \mbox{while $p(t)>0$} \,.
\]
Hence, $p'(t)<0$, while $p(t)>0$. So either $p(t)$ becomes zero at some $t_1$, or else $p(t)$ remains positive, and $\lim _{t \ra \infty} p(t)=b \geq 0$. 
In the latter case, $w(t)=w_0(t)+b+o(1)$, and then $t^{\alpha }a(t)=\frac{a_0+f(t)}{t^2}$, with $a_0=(2+\al)(n-2) \int _0^1 e^{(1-s)b} \, ds>0$, and $\lim _{t \ra \infty} f(t)=0$. Since $p(t)$ is bounded, 
by Lemma \ref{lma:5},  $p(t) \ra 0$ as $t \ra \infty$, i.e., $b=0$, and so $w(t)-w_0(t) \ra 0$, as $t \ra \infty$.
In case $p(t_1)=0$, we show by the same argument that the linear equation (\ref{36}) has either the second root at some $t_2>t_1$, or else $p(t)$ remains negative, and $\lim _{t \ra \infty} p(t)= 0$. (We have $\left(t^{n-1}p' \right)'>0$, when $p<0$, from which it is easy to deduce that $p(t)$ remains bounded.)
\medskip

Next, we  rule out the possibility of $p(t)$ keeping the same sign and tending to zero over an infinite interval $(t_k,\infty)$. We have
\[
t^{\al }a(t)=t^{\al }e^{w_0}\int_0^1  e^{(1-s) \left(w(t)- w_0(t) \right)} \, ds=\frac{(2+\al)(n-2)}{t^2} \left(1+o(1) \right), \, \mbox{as $t \ra \infty$} \,.
\]
Since Euler's equation (\ref{32}) has infinitely many roots on $(t_k,\infty)$, we conclude by Lemma \ref{lma:2} that $p(t)$ must vanish on that interval too. It follows that $p(t)$ changes sign infinitely many times. ($p(t)$ cannot remain positive,  so that it has its first root, after that $p(t)$ cannot remain negative,  so that it has its second root, and so on.)
\epf

We see that unlike the Gelfand equation ($\al=0$), infinitely many turns are possible for any $n \geq 3$, for large enough $\al$.
\medskip

We now turn to the case when $n \geq 10+4 \al$. In that case both roots of the characteristic equation of Euler's equation (\ref{32}) are negative, and so any solution of (\ref{32}) may have at most one root. By a simple comparison argument we shall show that $w(t)$ and $w_0(t)$ intersect at most twice, and then (in the case $\al =0$) we will show that $w(t)$ and $w_0(t)$ do not intersect at all. We have
\[
e^{w(t)}-e^{w_0(t)}> e^{w_0(t)} \left(w(t)-w_0(t) \right)=\frac{(2+\al)(n-2)}{t^{2+\al}}  \left(w(t)-w_0(t) \right) \,.
\]
Denoting $p(t)=w(t)-w_0(t)$, we then have from (\ref{31})
\beq
\lbl{40}
p''+\frac{n-1}{t}p'+ \frac{(2+\al)(n-2)}{t^2} p<0 \,.
\eeq
This inequality implies that $p(t)$ oscillates slower (faster) than $z(t)$, the solution  of Euler's equation (\ref{32}), provided that $p(t)<0$ ($>0$). The following lemma makes this observation precise.

\begin{lma}\lbl{lma:6}
Assume that $z(t)$ is a solution of (\ref{32}), such that $z(t_0)=p(t_0)$ and $z'(t_0)=p'(t_0)$ at some $t _0 \in (0,\infty)$, and $z(t)<0$ on $(t_0, \infty)$. Then $p(t)<0$ on $(t_0, \infty)$.
\end{lma}

\pf
Assume, on the contrary, that $p(\xi)=0$ at some $\xi \in (t_0, \infty)$, while $p(t)<0$ on $(t_0, \xi)$. From the equations (\ref{32}) and (\ref{40}) (keep in mind that $z(t)<0$)
\beq
\lbl{40a}
\left(pz'-p'z \right)'+\frac{n-1}{t} \left(pz'-p'z \right) <0 \s\mbox{on $(t_0, \xi)$} \,,
\eeq
and so the function $Q(t) \equiv t^{n-1} \left(pz'-p'z \right)$ is decreasing on $(t_0, \xi)$. But $Q(t_0)=0$, while $Q(\xi)=-\xi^{n-1} p'(\xi) z(\xi) \geq 0$, a contradiction.
\epf

\begin{thm}\lbl{thm:28} 
In case $n \geq 10+4 \al$ the solution curve of the generalized Gelfand's equation (\ref{30}) admits at most two turns.
\end{thm}

\pf
Let again $p(t)=w(t)-w_0(t)$. In view of Lemma \ref{lma:3}, we need to show that $p'(t)$ changes its sign at most twice, i.e., $p(t)$ changes its  monotonicity at most twice.
Since $p(t)$ satisfies the linear equation (\ref{36}), $p(t)$ cannot have points of positive local minimum, and of negative local maximum, and hence $p(t)$ changes its monotonicity once between two consecutive roots, and once  after its last root (since $p(t)$ tends to zero as $t \ra \infty$, which follows by Lemma \ref{lma:5}, the same way as in the proof of Theorem \ref{thm:1}), and no other changes of monotonicity are possible. We will show that $p(t)$ has at most two roots, which will imply that $p(t)$ changes its monotonicity at most twice.
\medskip

Assume that $p(t)$ has at least two roots, and let $t_1$ and $t_2$ denote the first two roots (if there are less than two roots, there are less than two turns). Then $p(t)$ is negative on $(0,t_1)$, positive on $(t_1,t_2)$, and again negative after $t_2$. Pick any point $t_0 \in (0,t_1)$, and let $Z(t)$ be the solution of Euler's equation (\ref{32}), such that $Z(t_0)=p(t_0)<0$, $Z'(t_0)=p'(t_0)$. We claim that $Z(t)$ vanishes on $(0,t_1)$. Indeed, assuming that $Z(t)<0$  on $(0,t_1)$, we argue as in Lemma \ref{lma:6}, and conclude that the function  $Q(t) \equiv t^{n-1} \left(pZ'-p'Z \right)$ is decreasing on $(t_0, t_1)$,  with  $Q(t_0)=0$, while $Q(t_1)=-t_1^{n-1} p'(t_1) Z(t_1) \geq 0$, a contradiction. At its root $Z(t)$ changes to being  positive, and it stays positive after its root, since $Z(t)$ is a solution of Euler's equation with two negative characteristic roots (the roots coincide when $n=10+4 \al$). In particular,
\beq
\lbl{40b}
Z(t)>0 \;\; \mbox{for $t>t_1$}\,.
\eeq

We now return to $p(t)$. After its second root $t_2$, it will either have the third root at some $t_3$, or stay negative and tend to zero as $t \ra \infty$. We now consider these cases in turn.
\medskip

\noindent
{\bf Case 1} $p(t_2)=p(t_3)=0$, $p(t)<0$ on $(t_2,t_3)$. In place (\ref{40a}), we now have (because of (\ref{40b}))
\beq
\lbl{40c}
\left(pZ'-p'Z \right)'+\frac{n-1}{t} \left(pZ'-p'Z \right) >0 \s\mbox{on $(t_1, \infty)$} \,.
\eeq
Integrating this over $(t_2,t_3)$, we get 
\[
-t_3^{n-1} p'(t_3) Z(t_3)+t_2^{n-1} p'(t_2) Z(t_2)>0 \,,
\]
which is a contradiction, since (using (\ref{40b})) both terms on the left are non-positive.
\medskip

\noindent
{\bf Case 2} $\;$ We have $p(t)>0$ on $(t_2,\infty)$, and $p(t) \ra 0$, as $t \ra \infty$. This case is possible (or rather, we are unable to rule this case out).  Then $p(t)$ changes its monotonicity twice, and the solution curve has two turns.
\epf

It is natural to expect that in the case $n \geq 10+4\al$ the solution curve does not turn at all (if $\al =0$, this is part of D.D. Joseph and T.S. Lundgren's result \cite{JL}). Surprisingly, we found this hard to prove. We shall prove this only if $\al =0$ (with a computer assistance at one step), so that we recover the classical result of   D.D. Joseph and T.S. Lundgren \cite{JL}). By a different method, the case $n \geq 10+4\alpha$ was covered by J. Jacobsen and K. Scmitt \cite{JS}, who proved that for $0<\lambda<(n-2)(\alpha +2)$, the problem (\ref{30}) has a unique solution, and no solution exists for $\lambda \geq (n-2)(\alpha +2)$. 
\medskip

So we  consider now the case $\al =0$. The linearized equation at $w_0(t)$ is then  Euler's equation (\ref{19}), whose characteristic exponents are given by (\ref{19.1}). In case $n \geq 10$, both characteristic exponents are negative and hence $z(t) \ra 0$ as $t \ra \infty$. This solution will either vanish once or keep the same sign depending on the initial conditions. Assume that the initial conditions are given at some $A>0$. By scaling we may assume that $z(A)=1$. The following lemma says that in order for $z(t)$ to vanish, $z'(A)$ must be negative and sufficiently large in absolute value.

\begin{lma}\lbl{lma:7}
For  $n \geq 10$ consider the problem
\beq
\lbl{41}
\s\s\s z''+\frac{n-1}{t}z'+  \frac{2n-4}{t^2}z=0 \,, \s z(A)=1, \s z'(A)=q \s\s ( \mbox{for $t>A>0$}) \,.
\eeq
Assume that 
\beq
\lbl{42}
qA>\frac{-n+2 - \sqrt{(n-2)(n-10)}}{2} \,.
\eeq
Then $z(t)>0$ on $(A,\infty)$.
\end{lma}

\pf
If $n=10$, the general solution of the equation in (\ref{41}) is
\[
z(t)=t^{-4} \left(c_1+c_2 \ln t \right) \,.
\]
From the initial conditions
\[
c_2=A^4(qA+4)>0 \,,
\]
since (\ref{42}) reads: $qA>-4$. If $c_1 \geq 0$, then $z(t)>0$ for all $t>0$. If $c_1 <0$, then the function $c_1+c_2 \ln t$ has a root, but it is smaller than $A$ (observe that $c_1+c_2 \ln A=A^4>0$). So that $z(t)>0$ on $(A,\infty)$.
\medskip
If $n>10$, the general solution of (\ref{41}) is
\[
z(t)=c_1 \left(\frac{t}{A} \right)^r+c_2 \left(\frac{t}{A} \right)^s=\left(\frac{t}{A} \right)^r \left(c_1+c_2 \left(\frac{t}{A} \right)^{s-r} \right) \,.
\]
where $r=\frac{-n+2 - \sqrt{(n-2)(n-10)}}{2}$ and $s=\frac{-n+2 + \sqrt{(n-2)(n-10)}}{2}>r$. From the initial conditions
\[
c_2=\frac{Aq-r}{s-r}>0 \,,
\]
in view of (\ref{42}). If $c_1 \geq 0$, then $z(t)>0$ for all $t>0$. If $c_1 <0$, then the function $c_1+c_2 \left(\frac{t}{A} \right)^{s-r}$ has a root, but it is smaller than $A$ (at $A$ this function equals $1$). So that $z(t)>0$ on $(A,\infty)$.
\epf

Let now $t_0$ denote the root of $w_0(t)$, i.e., $t_0=\sqrt{2n-4}$. With $p(t)=w(t)-w_0(t)$, we have $p(t_0)=w(t_0)<0$. 
By Lemma \ref{lma:6} we shall have $p(t)<0$ on $(t_0, \infty)$, provided the solution of the linearized equation (\ref{41}) with the initial conditions $z(t_0)=p(t_0)$ and $z'(t_0)=p'(t_0)$ satisfies $z(t)<0$ on $(t_0, \infty)$. By Lemma \ref{lma:7} (with $A=t_0$) this will happen if
\beq
\lbl{44}
\frac{t _0w'(t_0)+2}{w(t_0)}>
 \frac{-n+2 - \sqrt{(n-2)(n-10)}}{2} \,.
\eeq
When $n=10$, a numerical computation shows that the quantity on the left in (\ref{44}) is approximately $-1.72324$, while the one on the right is $-4$. When one increases $n$, numerical computations show that the quantity on the left in (\ref{44}) is monotone increasing, while the one on the right is decreasing rapidly. It follows that $p(t)<0$, i.e., 
 $w(t)<w_0(t)$ for all $t>0$.  We conclude that $p'(t)>0$ for $t > 0$, since we would get a contradiction in (\ref{35.1}) at any point of local maximum. From (\ref{17.1}), $\la '(t)>0$ for all $t>0$, i.e., the solution curve always travels to the right in the $(\la, u(0))$ plane.

\section{Micro-electromechanical systems (MEMS)}
\setcounter{equation}{0}
\setcounter{thm}{0}
\setcounter{lma}{0}

Recall that for the problem (with constants $\alpha \geq 0$, and $p>1$)
\beq
\lbl{20}
\s\s\s u''+\frac{n-1}{r}u' + \la \, \frac{r^{\alpha}}{(1-u)^p }=0 \s \mb{for $0<r<1$,} \s u'(0)=0, \s u(1)=0 \,
\eeq
the solution curve $(\la, u(0))$ is given by $\left(\frac{t ^{\al +2}}{w^{p+1}(t)}\,, 1-\frac{1}{w(t)} \right)$, parameterized by $t \in (0,\infty)$, where $w(t)$ is the solution of
\beq
\lbl{21}
w''+\frac{n-1}{t}w'=  \frac{t^{\al}}{w^p},   \s  \s w(0)=1, \s w'(0)=0 \s (t>0)\,,
\eeq
with $w(t)>0$ and $w'(t)>0$ for all $t>0$.
In particular, $\la (t)=\frac{t ^{\al +2}}{w^{p+1}(t)}$. 
Compute
\[
\la '(t)=\frac{t^{\al +1} w^p \left[ (\al+2)w-(p+1)tw' \right]}{w^{2(p+1)}} \,.
\]
To study the direction of the solution curve  we are interested in the sign of $\la '(t)$, or in the roots of the function $(\al+2)w(t)-(p+1)tw'(t)$. If we set this function to zero
\[
(\al+2)w-t(p+1)w'=0 \,,
\]
then the general solution of this equation is  $w(t)=ct^{\beta}$, with $\beta=\frac{\al+2}{p+1}$. One verifies that 
\[
w_0(t)=c_0 t^{\beta}, \s \mbox{with} \s c_0=\frac{1}{\left[\beta(\beta+n-2) \right]^{\frac{1}{p+1}}} \;.
\]
is a solution of the equation in (\ref{21}) (the guiding solution).
The linearized equation for (\ref{21}) is
\[
z''+\frac{n-1}{t}z'=-p t^{\al} w^{-p-1} z \,.
\]
At the solution $w=w_0(t)$, this becomes
\beq
\lbl{23}
z''+\frac{n-1}{t}z'+  \frac{p \beta(\beta+n-2)}{t^2}z=0 \,,
\eeq
which is again Euler's equation. Its characteristic equation has the roots
\[
r=\frac{-n+2 \pm \sqrt{(n-2)^2-4p \beta(\beta+n-2) }}{2} \,.
\]
The roots are complex, if $\beta$ satisfies
\beq
\lbl{24}
4p \beta ^2+4p(n-2) \beta -(n-2)^2>0 \,.
\eeq
When $n=1$, this happens when $\beta> \frac{p+\sqrt{p^2+p}}{2p}$, the larger root for the quadratic on the left, i.e., when $\al> \frac{(p+1)\left(p+\sqrt{p^2+p}\right)}{2p}-2$. If $p=2$, this becomes
 $\al>\frac{6+3\sqrt{6}}{4}-2$. This is the same inequality as $\al> -\frac{1}{2}+\frac{1}{2}\sqrt{\frac{27}{2}}$, obtained on p. 900 in J.A. Pelesko \cite{P}. When $n=2$, the inequality (\ref{24}) holds for all $\al \geq 0$, which will imply (as we show below) infinitely many turns of the solution curve, as was previously proved for the  case $p=2$ in \cite{P}, and for any $p>1$ in \cite{GW}.
\medskip

We now assume that $n \geq 3$. The inequality (\ref{24}) holds if $\beta$ is greater than the larger root for the quadratic on the left, i.e.,
for
\beq
\lbl{25}
\beta> \frac{n-2}{2} \, \frac{1}{p+\sqrt{p^2+p}} \,.
\eeq
This inequality is equivalent to the condition (\ref{0.5}).
\medskip

\noindent
{\bf Remark}
Since $\beta=\frac{\al+2}{p+1}>\frac{2}{p+1}$, we shall see that in the lower dimensions (\ref{25}) holds automatically, i.e., without restricting $\al$. Indeed, we rewrite the inequality
\[
\frac{2}{p+1}> \frac{n-2}{2} \, \frac{1}{p+\sqrt{p^2+p}}
\]
as
\beq
\lbl{26}
\frac{p+1}{p+\sqrt{p^2+p}}<\frac{4}{n-2} \,.
\eeq
On the left we have a decreasing function, which takes its maximum at $p=1$. So that (\ref{26}) will follow if
\[
\frac{2}{1+\sqrt{2}}<\frac{4}{n-2} \,,
\]
which happens for $3 \leq n \leq 6$. 
\medskip

Similarly to Lemma \ref{lma:3} we prove the following lemma.
\begin{lma}\lbl{lma:4}
Assume that $w(t)$ and $w_0(t)$ intersect infinitely many times. Then the solution curve of (\ref{20}) makes infinitely many turns.
\end{lma}

\begin{thm}
Assume that $\beta=\frac{\al+2}{p+1}$ satisfies (\ref{24}) (i.e., either $n=2$, or $n \geq 3$, and the condition (\ref{0.5}) holds).
Then the solution curve of (\ref{20}) makes infinitely many turns.
\end{thm}

\pf
In view of the Lemma \ref{lma:4}, we need to show that $w(t)$ oscillates infinitely many times around $w_0(t)$. We show first  that these functions get close to each other, as $t$ increases. Denote $P(t)=w(t)-w_0(t)$. It satisfies
\beq
\lbl{26a}
P''+\frac{n-1}{t}P'+a(t) P=0 \,,
\eeq
where $a(t)=p \, t^{\al } \int_0^1 \frac{1}{\left[s w(t)+(1-s) w_0(t) \right]^{p+1}} \, ds>0$. As in the proof of the Theorem \ref{thm:1}, we see that either $P(t)$ has infinitely many roots, or else $P(t)$ keeps the same sign over some infinite interval $(t_k, \infty)$, and tends to a constant as $t \ra \infty$. We now rule out the latter possibility.  Write
\[
a(t)=p \, t^{\al }\frac{1}{w_0^{p+1}}\int_0^1  \frac{1}{ \left[s \frac{w(t)}{w_0(t)} +(1-s)   \right]^{p+1}} \, ds=\frac{p \beta(\beta+n-2)}{t^2} \left(1+o(1) \right)  \,,
\]
as $t \ra \infty$. (Observe that $\frac{w(t)}{w_0(t)}=1+\frac{P(t)}{w_0(t)} \ra 1$, as $t \ra \infty$.)
Since Euler's equation (\ref{23}) has infinitely many roots on $(t_k,\infty)$, we conclude by Lemma \ref{lma:2} that $P(t)$ must vanish on that interval too. It follows that $P(t)$ changes sign infinitely many times.
\epf

\begin{thm} 
Assume that $n \geq 2+ \frac{2(\al+2)}{p+1} \, \left(p+\sqrt{p^2+p} \right)$. Then  the solution curve of  (\ref{20}) admits at most two turns.
\end{thm}

\pf
We follow the proof of the Theorem \ref{thm:28}.
With $P(t)=w(t)-w_0(t)$, we need to show that $P(t)$ changes its monotonicity at most twice.
Since $P(t)$ satisfies the linear equation (\ref{26a}), $P(t)$ cannot have points of positive local minimum, and of negative local maximum, and hence $P(t)$ changes its monotonicity once between two consecutive roots, and once  after its last root (since $z(t)$ tends to zero as $t \ra \infty$), and no other changes of monotonicity are possible. We will show that $P(t)$ has at most two roots, which will imply that $P(t)$ changes its monotonicity at most once. Under our conditions  Euler's equation (\ref{23}) has at most one root, while $P(t)$ satisfies 
\[
P''+\frac{n-1}{t}P'+  \frac{p \beta(\beta+n-2)}{t^2}P<0 \,,
\]
since the nonlinearity in (\ref{20}) is convex in $u$. The rest of the proof is identical to that of the Theorem \ref{thm:28}.
\epf

\section{Morse index of solutions to the Gelfand and MEMS problems}
\setcounter{equation}{0}

We now use the generating solution to show that all turning points of the Gelfand problem are non-degenerate, and that the Morse index of solutions increases by one at each turning point, thus recovering a result of K. Nagasaki and T.  Suzuki \cite{N2}. Recall that positive solutions of the Gelfand problem
\beq
\lbl{60}
\Delta u+\la e^u=0, \s \mbox{for $|x|<1$}, \s u=0 , \s \mbox{when $|x|=1$} 
\eeq
are radially symmetric, i.e., $u=u(r)$, $r=|x|$, with $u'(r)<0$, and that the solution curve in the $(\la ,u(0))$ plane is given by (\ref{0.3}) and (\ref{0.21}). In particular, $\la (t)=t^2e^{w(t)}$, where $w(t)$ is the solution of (\ref{0.21}), the generating solution.

\begin{thm}\lbl{thm:50}
Let $u(t_n)$ be a singular solution of (\ref{60}), i.e., $\la '(t_n)=0$. Then $u(t_n)$ is non-degenerate, i.e., $\la ''(t_n) \ne 0$.
\end{thm}

\pf
Compute $\la '(t)=te^{w(t)} \left(2+t w'(t) \right)$. Since $\la '(t_n)=0$, we have $2+t_n w'(t_n)=0$. Then
\beq
\lbl{61}
\la ''(t_n)=t_ne^{w(t_n)} \left( w'(t_n)+t_n w''(t_n)\right),
\eeq
and we need to show that $w'(t_n)+t_n w''(t_n) \ne 0$. For the guiding solution $w_0(t)=\ln(2n-4) -2 \ln t$ we have
\[
w_0'(t_n)+t_n w_0''(t_n) =0 \,.
\]
Since $w'(t_n)=w_0'(t_n)$ ($=-\frac{2}{t_n}$), it suffices to show that
\beq
\lbl{62}
w''(t_n) \ne w_0''(t_n) \,.
\eeq
The function $p(t)=w_0(t)-w(t)$ satisfies the linear equation (\ref{36}), with $p'(t_n)=0$. It follows that $w(t_n) \ne w_0(t_n)$, since otherwise $p(t) \equiv 0$, which is impossible. Since both $w(t)$ and $w_0(t)$ satisfy the same equation (\ref{0.21}), we conclude (\ref{62}) by expressing the second derivatives from the corresponding equations.
\epf

By C.S. Lin and W.-M. Ni \cite{LN}, any solution of the linearized problem for (\ref{60}) is radially symmetric, and hence it satisfies
\beq
\lbl{63}
\omega ''+\frac{n-1}{r}\omega'+\la e^u \omega=0, \s \mbox{for $0<r<1$}, \s \omega'(0)=\omega(1)=0 \,.
\eeq
We call $u(r)$ a singular solution of (\ref{60}) if the problem (\ref{63}) has a non-trivial solution. (Differentiating (\ref{60}) in $t$, and setting $t=t_n$, it is easy to see that  a  solution is singular  iff $\la '(t_n)=0$.) The following lemma was proved in P. Korman \cite{K}.
\begin{lma}\lbl{lma:60}
Let $u(r)$ be a singular solution of (\ref{60}). Then 
\[
\om (r) =ru'(r)+2
\]
gives a solution of (\ref{63}).
\end{lma}

We now recover the following  result of K. Nagasaki and T.  Suzuki \cite{N2}.

\begin{thm}\lbl{thm:60}
As one follows the solution curve of (\ref{60}) in the direction of increasing $u(0)$, the Morse index of solution increases by one at each turn.
\end{thm}

\pf
The Morse index of solution is the number of negative eigenvalues $\mu$ of 
\[
\Delta \om+\la e^u \om+\mu \om =0, \s \mbox{for $|x|<1$}, \s \om=0 , \s \mbox{when $|x|=1$}  \,.
\]
By \cite{LN} solutions of this problem are radially symmetric. At a singular solution $\mu=0$, and then $\om (r) =ru'(r)+2$ by Lemma \ref{lma:60}.
Assume that at a singular solution $u(t_n)$, $\mu (t_n)=0$ is the $k$-th eigenvalue. Following \cite{N2}, we will show that $\mu'(t_n)<0$, which means that for $t<t_n$ ($t>t_n$) the $k$-th eigenvalue is positive (negative), i.e., the Morse index increases by one through $t=t_n$. We shall show that the sign of $\mu'(t_n)$ is the same as that of $-\left(\la ''(t_n) \right)^2$, which is negative by the Theorem \ref{thm:50}. Recall the following known formulas (here $u=u(t_n)$, $\om$ is a solution of (\ref{63}), and $B$ is the unit ball around the origin in $R^n$):
\[
\mu'(t_n) \int _B \om ^2 \, dx=-\la (t_n) \int _B e^u \om ^3 \, dx \s (\mbox{p. $11$  in \cite{K1}}) \,,
\]
\[
-\la (t_n) \int _B e^u \om ^3 \, dx=\la ''(t_n) \int _B f(u) \om \, dx \s (\mbox{p. $3$   in \cite{K1}}) \,,
\]
\[
\int _B f(u) \om \, dx=\frac{1}{2 \la (t_n)} u'(1) \om '(1) \s (\mbox{p. 5  in \cite{K1}}) \,.
\]
Since $u'(1)<0$, it follows from these formulas that the sign of $\mu'(t_n)$ is opposite to that of $\la ''(t_n) \, \om '(1)$. Using the Lemma \ref{lma:60}, we have $\om'(r)=u'+ru''$, $\om'(1)=u'(1)+u''(1)$. Recall that $u(r)=w(t)+a$, with $t=br$ (where $w(t)$ is the generating solution). Observing that $b=t_n$ for $r=1$, we have $u(r)=u(0)+w(t_n r)$, and then 
\[
\om '(1)=t_n \left(w'(t_n)+t_n w''(t_n) \right),
\]
which by (\ref{61}) has the same sign as $\la ''(t_n)$. It follows that the sign of $\mu'(t_n)$ is the same as that of $-\left(\la ''(t_n) \right)^2<0$.
\epf

For the MEMS problem ($p>0$)
\beq
\lbl{70}
\Delta u+\la \frac{1}{(1-u)^p }=0, \s \mbox{for $|x|<1$}, \s u=0 , \s \mbox{when $|x|=1$}
\eeq
the situation is similar, which is a new result.

\begin{thm}\lbl{thm:70}
Let $u(t_n)$ be a singular solution of (\ref{70}), i.e., $\la '(t_n)=0$. Then $u(t_n)$ is non-degenerate, i.e., $\la ''(t_n) \ne 0$. Moreover, when  one follows the solution curve of (\ref{70}) in the direction of increasing $u(0)$, the Morse index of solution increases by one at each turn.
\end{thm}

As before, by \cite{LN}  solutions of the linearized problem corresponding to  (\ref{70}) are radially symmetric, and hence they satisfy
\beq
\lbl{71}
\s\s \omega ''+\frac{n-1}{r}\omega'+\la \frac{p}{(1-u)^{p+1} } \omega=0, \s \mbox{ $0<r<1$}, \s \omega'(0)=\omega(1)=0 \,.
\eeq

\begin{lma}\lbl{lma:70}
Let $u(r)$ be a singular solution of (\ref{70}). Then 
\[
\om (r) =ru'(r)-\frac{2}{p+1} u(r)+\frac{2}{p+1}
\]
gives a solution of (\ref{71}).
\end{lma}

\pf
The function $v(r) \equiv ru'(r)-\frac{2}{p+1} u(r)+\frac{2}{p+1}$ solves (\ref{71}), and we have $v'(0)=0$, $v(0)>0$ (since solutions of (\ref{70}) are smaller than $1$). By scaling of $\omega (r)$, we may assume that $\omega (0)=v(0)$, and then by the uniqueness result for this type of initial value problems (see \cite{PS1}), it follows that $\omega (r) \equiv v(r)$.
\epf

\noindent
{\bf Proof of the Theorem \ref{thm:70}} The proof is similar to that of the Theorem \ref{thm:60}. This time $\la (t)=t^2w^{-p-1}(t)$.
Compute 
\[
\la '(t)=tw^{-p-2}(t) \left(2w(t)-(p+1)t w'(t) \right). 
\]
Since $\la '(t_n)=0$, we have 
\beq
\lbl{71.1}
2w(t _n)-(p+1)t _n w'(t _n)=0 \,. 
\eeq
Then
\beq
\lbl{72}
\la ''(t_n)=-t_nw^{-p-2}(t_n) \left( (p-1)w'(t_n)+(p+1)t_n w''(t_n) \right),
\eeq
and we need to show that 
\[
S \equiv (p-1)w'(t_n)+(p+1)t_n w''(t_n) \ne 0 \,,
\]
to conclude  that $\la ''(t_n) \ne 0$. 
Using the equation (\ref{21}), and then (\ref{71.1}), we express
\[
S=-\left[(p+1)(n-1)-(p-1) \right]w'(t_n)+\frac{(p+1)t_n}{w^p(t_n)}
\]
\[
=-\frac{2 \left[(p+1)(n-1)-(p-1) \right]w(t_n)}{(p+1)t_n}+\frac{(p+1)t_n}{w^p(t_n)} \,.
\]
For the guiding solution $w_0(t)=c_0t^{\beta}$, $\beta=\frac{2}{p+1}$, we have 
\[
-\frac{2 \left[(p+1)(n-1)-(p-1) \right]w_0(t_n)}{(p+1)t_n}+\frac{(p+1)t_n}{w_0^p(t_n)}=0 \,.
\]
Observing that $S$ is a decreasing function of $w(t_n)$
 we conclude that $S \ne 0$, once we  show that
\beq
\lbl{620}
w(t_n) \ne w_0(t_n) \,.
\eeq
The function $p(t)=w_0(t)-w(t)$ satisfies the linear equation (\ref{26a}), with $p'(t_n)=0$. Then (\ref{620}) is true, since otherwise $p(t) \equiv 0$, which is impossible. 
\medskip

We see as before that the sign of $\mu'(t_n)$ is opposite to that of $\la ''(t_n) \, \om '(1)$. By Lemma \ref{lma:70}, 
\[
\omega'(1)=\frac{p-1}{p+1} u'(1)+u''(1) \,.
\]
In terms of the generating solution $w(t)$, we have $u(r)=1-a w(t)$, $t=br$, with $b=t_n$ at the singular solutions. Then $u(r)=1-aw(t_n r)$, and we have 
\[
 \omega'(1)
 =-\frac{at_n}{p+1} \left[ (p-1) w'(t_n)+(p+1)t_nw''(t_n) \right]\,. 
\]
Comparing this with (\ref{72}), we conclude that  $\omega'(1)$ has the same sign as $\la ''(t_n)$, and then the sign of $\mu'(t_n)$ is the same as that of $-\left(\la ''(t_n) \right)^2<0$, concluding the proof as before.
\epf

\section{The Henon equation}
\setcounter{equation}{0}
\setcounter{thm}{0}
\setcounter{lma}{0}

We study positive solutions of the Dirichlet problem for the  Henon equation
\beq
\lbl{h1}
u''+\la |x|^{\al} u^p=0, \s -1<x<1, \s u(-1)=u(1)=0 \,.
\eeq
Here $\al >0$ and $p>1$ are constants, $\la>0$ is a parameter. This problem has both symmetric (even) and non-symmetric positive solutions, see D. Smets et al \cite{S}, R. Kajikiya \cite{J}.
The exact multiplicity of the positive solutions is not known. We shall approach this problem by using ``shooting", scaling and numerical computations.
\medskip

For a variable $\xi >0$ we consider the initial value problem
\beq
\lbl{h2}
z''+ |x|^{\al} z^p=0, \s z(\xi)=1, \s  z'(\xi)=0 \,.
\eeq
Let $a(\xi)$ denote the first root of $z(x)$ which is greater than $\xi$, and $b(\xi)$ the first root of $z(x)$ which is to the left of $\xi$.

\begin{thm}
Assume that the equation
\beq
\lbl{h3}
a(\xi)=-b(\xi) \s\s ( \xi >0)
\eeq
has a unique solution $\xi _0 >0$. Then for any positive $\la$ the problem (\ref{h1}) has exactly three positive solutions: $u_1(x)$ which is an even function, $u_2(x)$ which has its point of maximum at $\xi=\frac{\xi _0}{a(\xi _0)}$, and $u_3(x)=u_2(-x)$. Moreover, if $b$ denotes the maximum value of $u_2(x)$ (i.e., $b=u_2(\xi)$), then $\la=\frac{a(\xi _0)^{\al +2}}{b^{p-1}}$.
\end{thm}

\pf
Denote $\eta =a(\xi _0)$. In (\ref{h2}) we let $x=\eta t$, $z=\frac{1}{b} v$, obtaining
\[
v''+ \frac{\eta^{\al +2}}{b^{p-1}} |t|^{\al} v^p=0=0, \s v(\pm 1)=0 \,,
\]
i.e., $v(t)$ is a solution of (\ref{h1}), with $\la = \frac{\eta^{\al +2}}{b^{p-1}}$. The maximum value of this solution is equal to $b$, and it occurs at $\xi= \frac{\xi _0}{\eta}$. Solutions of (\ref{h1}) at other $\la$'s are obtained by scaling of $u$. 
\medskip

We now show that there is exactly one positive solution, which takes its maximum value at a positive $x$. Let $u_3(x)$ be another solution of (\ref{h1}), with the maximum value achieved at $\xi _1 >0$. Assume first that $\xi _1 =\xi _0$. By scaling of $u$ and $x$, we obtain from $u_3(x)$ a solution of (\ref{h1}) (at the same $\la$), which at $\xi _0$ has the same initial data as $u_1(x)$ (and so is identical to $u_1(x)$), but it has its first root at some $x \ne 1$, a contradiction. Next, assume that $\xi _1 \ne \xi _0$. Again, we scale $u=Av$, $x=Bt$. Choose $A=u(\xi _1)$, then $v(\xi _1)=1$. Now choose $B$, so that $\la B^{2+\al}A^{p-1}=1$. Then we get a solution of (\ref{h2}), for which $a(\xi _1)=\frac1B$, $b(\xi _1)=-\frac1B$, a contradiction with the uniqueness of solution of (\ref{h3}).
\epf

For particular $\al$ and $p$ one can verify computationally that (\ref{h3}) has a unique solution.
\medskip

\noindent
{\bf Example} $\al=2$, $p=3$. The graphs of $a(\xi)$ and $-b(\xi)$, computed using {\em Mathematica}, are given in the Figure $2$. For $\xi>0$ and small, we have  $-b(\xi)>a(\xi)$, as a graph on a smaller scale shows. We see that the graphs of $a(\xi)$ and $-b(\xi)$ intersect exactly once for $\xi >0$.
This computation provides a computer assisted proof that the problem
\[
u''+\la x^{2} u^3=0, \s -1<x<1, \s u(-1)=u(1)=0
\]
has exactly three positive solutions for any $\la > 0$.

\begin{figure}
\begin{center}
\scalebox{0.9}{\includegraphics{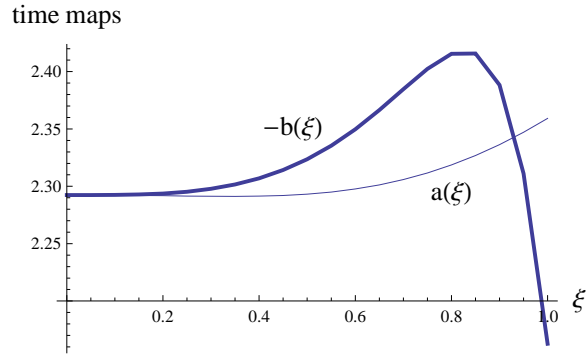}}
\caption{  The functions $a(\xi)$ and $-b(\xi)$}
\end{center}
\end{figure}

\medskip

We have seen a similar picture for many other $\al >0$ and $p>1$ that we tried. For particular $\al $ and $p$ one can produce a computer assisted proof, based on these calculations. However, some restrictions on $\al $ and $p$ appear to be necessary. For example for $\al=p=2$ the graph of $-b(\xi)$ is below that of $a(\xi)$, which indicates that there are no symmetry breaking solutions. In a recent paper S. Tanaka \cite{T} proved that symmetry breaking solutions exist, provided that $\al (p-1) \geq 4$. For symmetric (even) solutions existence and uniqueness is known for all $\al >0$ and $p>1$, see e.g., \cite{K1}.
\medskip

\noindent{\bf Acknowledgments} $\;$ It is a pleasure to thank Ken Meyer for many stimulating discussions. I wish to thank the referee for a number of useful comments, and for catching an error in an earlier version of the proof of Theorem 5.3.

%\begin{figure}
%\scalebox{0.9}{\includegraphics{}}
%\caption{  Solution curve for the problem (\ref{58})}
%\end{figure}

\end{document}